\documentclass[10pt]{amsart}
\usepackage{amssymb,amsmath}
\usepackage[all]{xy}
\usepackage{graphicx}
\numberwithin{equation}{section}

\newtheorem{theorem}{Theorem}[section]
\newtheorem{definition}[theorem]{Definition}
\newtheorem{corollary}[theorem]{Corollary}

\newtheorem{lemma}[theorem]{Lemma}
\newtheorem{remark}[theorem]{Remark}

\newcommand{\RR}{\mathbb{R}}

\newcommand{\kk}{\boldsymbol{k}}
\newcommand{\ZZ}{\mathbb{Z}}

\newcommand{\NN}{\mathbb{N}}
\newcommand{\AI}{A_\infty}

\newcommand{\bx}{\boldsymbol{x}}

\newcommand{\WH}[1]{\widehat{#1}}
\newcommand{\UL}[1]{\underline{#1}}
\newcommand{\NOV}{\Lambda_{nov}}
\begin{document}
\title{Strong homotopy inner product of an $\AI$-algebra }
\author{Cheol-Hyun Cho}

\address{Department of Mathematical Sciences, Seoul National University,
San 56-1, Shinrimdong, Gwanak-gu, Seoul, South Korea,  Email address : chocheol@snu.ac.kr }

\begin{abstract}
We introduce a strong homotopy notion of a cyclic symmetric inner product of an $\AI$-algebra
and prove a characterization theorem in the formalism of the infinity inner products by Tradler.
We also show that it is equivalent to the notion of a non-constant symplectic structure on the corresponding
formal non-commutative supermanifold.

We show that (open Gromov-Witten type) potential for a cyclic filtered $\AI$-algebra is invariant under the
cyclic filtered $\AI$-homomorphism up to reparametrization, cyclization and a constant addition, generalizing the work of Kajiura. 
\end{abstract}
\maketitle

\bigskip

\section{Introduction}
Inner product is an important notion in general and the case of
a strong homotopy associative algebra ($\AI$-algebra) is not an exception. 
 A suitable notion of an inner product in this case
is the cyclic symmetric inner product introduced by Kontsevich \cite{K}, where the $\AI$-operations, say $m_k:A^{\otimes k} \to A$, and
a cyclic symmetric inner product is required to satisfy the following identities: For $x_i \in A$,
\begin{equation}
	<m_k(x_1,\cdots,x_k),x_{k+1}> = \pm <m_k(x_2,\cdots,x_{k+1}),x_{1}>.
\end{equation}
The exact sign will be explained in the next section due to different conventions used in literatures.
This notion for the $\AI$-algebras and $\AI$-categories is crucial, for example, in the
work of Kontsevich-Soibelman\cite{KS} or  Costello\cite{Cos}. 

Another important application is to define a potential for an $\AI$-algebra.
Unlike the closed Gromov-Witten invariants,
the general open (genus zero) Gromov-Witten numbers counting pseudo-holomorphic discs with
intersection  conditions, does depend on several choices and does not define an invariant in general.
(But for the Lagrangian submanifold given by the
fixed points of an anti-holomorphic involution with certain conditions, Welshinger has defined
an invariant count. See \cite{W},\cite{C},\cite{S}).

In general, one should consider a filtered $\AI$-algebra of a Lagrangian submanifold, where
different choices made during the construction give rise to the same $\AI$-algebra up to $\AI$-homotopy equivalences.
But to consider the relation of their structure constants, one need a finer structure, namely cyclic filtered $\AI$-algebras
and cyclic filtered $\AI$-homomorphisms. 

For a cyclic $\AI$-algebra, say $A$, one can define a potential function (Definition \ref{podef}) such that
for any cyclic $\AI$-algebra $B$ which is cyclic $\AI$-homomorphic to $A$, their respective potential functions are related by reparametrizaion
and cyclization, by the work of Kajiura \cite{Kaj}. We define corresponding notions in the filtered case and extend his results to the case of a
filtered $\AI$-algebra (Theorem \ref{thm:pofil}), which is a natural setting for the open Gromov-Witten potential.
We remark that one should consider {\it cyclic } filtered $\AI$-algebras for the purpose of counting $J$-holomorphic discs since
$\AI$-algebra which is not cyclic, is not enumerative. For a cyclic filtered $\AI$-algebra, even though
the individual constants do not define invariants in general, the potential may be considered to be invariant
(up to reparametrization and up to a constant addition), since a cyclic filtered $\AI$-homomorphism
preserves the potential in this sense.

So far, we have addressed the importance of the cyclicity for the $\AI$-structures.
The main part of this paper concerns the homotopy notions of a cyclic symmetric inner product of an $\AI$-algebra.
Tradler\cite{T} has defined the notion of an infinity inner product as an $\AI$-bimodule map from an $\AI$-algebra $A$
to its dual $A^*$. Cyclic symmetric inner products also can be written in this formalism (Lemma \ref{lemma:cyc}).

 The main question we address in this paper is with what condition does an $\AI$-algebra $A$ has a homotopy equivalent $\AI$-algebra $B$ such that
 $B$ has a cyclic symmetric inner product. Roughly speaking, we say that such an $\AI$-algebra has a {\em strong homotopy inner product}.
(Definition \ref{shidef}). We show that the existence of a strong homotopy inner product is equivalent to the following conditions on an $\AI$-bimodule map $\phi:A \to A^*$: skew symmetry , closedness condition and homological non-degeneracy. This is proved in the characterization Theorem \ref{thm:sh}.
The skew symmetry property only concerns the $\AI$-bimodule structures (Theorem \ref{thm:hi}). This is good enough
in the case of Tradler and Zeinalian to associate two Hochschild cohomologies $H^*(A,A)$ and $H^*(A,A^*)$
(to obtain a BV structure in the former)  where
the relation really concerns the $\AI$-bimodule structure of $A$ and $A^*$.
 
The strong homotopy inner product actually turns out to be equivalent to the notion of  the non-constant symplectic structure
on the corresponding formal non-commutative supermanifold(Theorem \ref{thm:equiv}). And the characterization theorem \ref{thm:sh} proves 
the equivalence of these two notions.

The question we hope to answer in the future is whether the construction of a filtered $\AI$-algebra of Lagrangian submanifolds by 
Fukaya, Oh, Ohta and Ono can be made cyclic symmetric. Due to the heavy use of abstract perturbations, 
the current construction so far is not cyclic symmetric, although there was an announcement by Fukaya that
in a different setting of De Rham chains, it can be made cyclic symmetric. If one shows that these two constructions
are $\AI$-homotopy equivalent ( which was proved in the case without quantum contribution in \cite{FOOO}), 
then, this would provide a strong homotopy inner product on the $\AI$-algebra on singular chains constructed in \cite{FOOO}.

We would like to thank Hiroshige Kajiura for very helpful conversations, and for the comments on the draft.

\section{Sign conventions and grading}
First, we explain the sign conventions regarding $\AI$-algebras and cyclic symmetry.
There are usually two conventions regarding $\AI$-formulas, where two sign conventions differ due to the degree shifting of
the vector spaces.
Before the degree shift,
\begin{definition}\label{def1}
 An $\AI$-algebra  $(A,\{m^{ns}_*\})$ (with no shifting) consists of a $\ZZ$-graded vector space $A$ with
a collection of multilinear maps $m:= \{ m^{ns}_n:A^{\otimes n} \to A\}_{n\geq 1}$ of
degree $2-n$ satisfying the following equation for each $k=1,2,\cdots.$
\begin{equation}
0=\sum_{k_1+k_2=k+1} \sum_{i=1}^{k_1-1} (-1)^{\epsilon_1} m^{ns}_{k_1}(x_1,\cdots,x_{i-1},m_{k_2}^{ns}(x_i,\cdots,x_{i+k_2-1}),\cdots,x_k)
\end{equation}
where $\epsilon_1= i(k_2+1)+ k_2(|x_1|+ \cdots + |x_{i-1}|)$.
Here $ns$ in the superscript means `no shifting'.
\end{definition}
\begin{remark}
Here, we only consider a strict $\AI$-algebra over the field $\kk$ with $char(\kk)=0$.
In particular,  we assume $m_0=0$. The case of $m_0 \neq 0$ is discussed in the last section of this paper.
\end{remark}

Instead, we will use the following definition of $\AI$-algebra after the degree shifting.
\begin{definition}\label{def2}
 An $\AI$-algebra $(A,\{m_*\})$ consists of a $\ZZ$-graded vector space $A$ with
a collection of multilinear maps  $m:= \{ m_n:A[1]^{\otimes n} \to A[1]\}_{n\geq0}$ of
degree one satisfying the following equation for each $k=1,2,\cdots.$
\begin{equation}
0=\sum_{k_1+k_2=k+1} \sum_{i=1}^{k_1-1} (-1)^{\epsilon_1} m_{k_1}(x_1,\cdots,x_{i-1},m_{k_2}(x_i,\cdots,x_{i+k_2-1}),\cdots,x_k)
\end{equation}
where $\epsilon_1= |x_1|'+ \cdots + |x_{i-1}|'$.
\end{definition}
Here $|x_i|$ is the degree of the element $x_i$, and $|x_i|'$ is the shifted degree.
Hence $|x_i|=|x_i|'+1$.

The above two conventions are related by the formula
\begin{equation}
m_{k}^{ns}(x_1,\cdots,x_k) = (-1)^{\sum_{i=1}^{k-1} (k-i)|x_i|} m_{k}(x_1,\cdots,x_k).
\end{equation}
The first convention (Definition \ref{def1}) has an advantage that if $m_k \equiv 0$ for $k \geq 3$, then
$A$ has the structure of a differential graded algebra with a differential $m_1$ and a product $m_2$
(i.e. signs match).
For the second sign convention (Definition \ref{def2}), this is not true since signs do not match, 
but this convention has an advantage that the sign rules in this case are just the Koszul sign convention. In this paper, we will
follow the second convention, which will simplify sign considerations.
We refer readers to \cite{GJ},\cite{T} for details regarding sign convention.

The signs for an inner product are also changed as follows.
First consider the case that inner product is defined by $m_2$.  
If the product $m^{ns}_2$ (or $m^{ns}_{2,0}$ in the filtered case) is graded symmetric, then
it turns out that $m_2$ (or $m_{2,0}$) is skew symmetric, which can be easily checked
from the relation $m_2^{ns}(a,b) = (-1)^{|a|}m_{2}(a,b)$.
\begin{lemma}
If $m^{ns}_2$ (or $m^{ns}_{2,0}$) is graded symmetric, then $m_2$ (or $m_{2,0}$) is skew symmetric. Namely, we have
$$m_2(a,b) + (-1)^{|a|'|b|'}m_2(b,a)=0.$$
\end{lemma}
In general, we can set $$<a,b>^{ns}=(-1)^{|a|}<a,b>$$
Then, symmetry of $<,>$ on a vector space corresponds to the skew symmetry on the shifted vector space.

Now we recall the notion of a cyclic symmetric $\AI$-algebra, which was first introduced by Kontsevich
as a non-commutative analogue of a symplectic structure (\cite{K}, also see section \ref{sec:nonc}).
\begin{definition}
A strict $\AI$-algebra $(A,\{m_*\})$ is said to have a {\it cyclic symmetric} inner product if
there exists a skew symmetric non-degenerate, bilinear map $$<,> : A[1] \otimes A[1] \to \kk$$
such that for all integer $k \geq 1$,
\begin{equation}
	<m_k(x_1,\cdots,x_k),x_{k+1}> = (-1)^{K(\vec{x})}<m_k(x_2,\cdots,x_{k+1}),x_{1}>.
\end{equation}
Here, $(-1)^{K(\vec{x})}$ denotes the sign given by Koszul sign convention.
Namely,
\begin{equation}
(-1)^{K(\vec{x})} = (-1)^{|x_1|'(|x_2|' + \cdots +|x_{k+1}|')}.
\end{equation}
For short, we will call such an algebra, cyclic $\AI$-algebra.
\end{definition}
\begin{remark}
This may be considered as a generalization of a property of a Frobenious algebra. 
Let $B$ be a Frobenious algebra (commutative) with the inner product $<a,b> = \theta(ab)$
given by a functional $\theta$ on $B$. Then, we have
$$<ab,c>= <a,bc>= <bc,a> = \theta(abc). $$
\end{remark}

Without the degree shifting, the cyclic symmetry equation becomes
\begin{equation}
	<m^{ns}_k(x_1,\cdots,x_k),x_{k+1}> = (-1)^{k+ K(\vec{x})}<m^{ns}_k(x_2,\cdots,x_{k+1}),x_{1}>.
\end{equation}
Note the additional sign depending on $k$. But $(-1)^k$ cannot be avoided. For example, when $k=1$, one can easily check
(from the Leibniz rule) that singular cohomology with Poincare pairing satisfies
$$ < m^{ns}_1 x, y> = - (-1)^{|x||y|} < m^{ns}_1 y, x>$$

For the signs in the equations, we will write $K$ throughout the paper whenever the Koszul convention is used.

We say that an inner product has degree $\alpha$ if $<a,b> \neq 0$ implies $deg(a) + deg(b) + \alpha =0$.
Inner product coming from geometry has degree $-n$ which is the dimension of the manifold.
As an inner product in a shifted vector space $A[1]$, it has degree $(-n+2)$.

%

\section{$\AI$-bimodule and infinity inner products}
In this section, we first recall the Tradler's definition of an infinity inner product of an $\AI$-algebra, and we define
the notions of (strong) homotopy inner products of an $\AI$-algebra in this setting.

To understand the Tradler's definition, we first consider the case of a commutative Frobenious algebra.
Let $(A, \cdot)$ be an algebra over the field $\kk$, let $A^*=Hom(A,\kk)$ be its dual.
Then, $A^*$ has an $A$-$A$-bimodule structure given as follows.
Namely, for $v^* \in A^*, a_1,a_2, w\in A$, we define $a_1 \cdot v^* \cdot a_2 \in A^*$ such that
\begin{equation}
	\big( a_1  \cdot v^* \cdot a_2 \big) (w) = v^*(a_2 \cdot w \cdot a_1).
\end{equation}
The algebra $A$ has a $A$-$A$-bimodule structure clearly, and consider the
$A$-$A$-bimodule homomorphism $\phi:A \to A^*$, and define $<a,b> = \phi(a)(b)$.
We also assume that $<a,b>=<b,a>$.
The map $\phi$ being a bimodule map implies that
\begin{eqnarray*}
<a \cdot b, c> &=& \phi(a \cdot b)(c) \\
&=& (a \cdot \phi ) (b)(c) \\
&=& \phi(b) (c \cdot a) = <b, c\cdot a>
\end{eqnarray*}

Similarly, one has the identity $<b \cdot a, c> = <b,a \cdot c>$. Hence the Frobenious property
may be written in terms of maps $\phi : A \to A^*$. Generalizing this idea, Tradler \cite{T} has defined the notion of an infinity inner product
of an $\AI$-algebra as an $\AI$-bimodule map $\psi :A \to A^*$ where $A^*$ is an $\AI$-bimodule over $A$.
Although infinity inner product is a generalized notion of an inner product, its relationship
with cyclic symmetric inner product has not been known. In this paper, we define homotopy notions
of cyclic symmetric inner products and prove a classification theorem in terms of the Tradler's definition.
It turns out that this provides a direct link to the non-constant symplectic structure (which is
a homotopy notion of a constant symplectic structure via a non-commutative Darboux theorem) which was introduced by
Kajiura.

First, we write cyclic symmetry in this formalism, and define a cyclic $\AI$-homomorphism, and
define homotopy notions of a cyclic symmetric inner product.

There exists a canonical $\AI$-bimodule structure on the dual $A^*$ of an $\AI$-algebra $(A,m)$.
Recall that $\AI$-bimodule structure of $M$ over $A$ is described by sequence of maps of degree one
$$b_{k,l}:A[1]^{\otimes k} \otimes M \otimes A[1]^{\otimes l} \to M,$$
satisfying $\WH{b} \circ \WH{b}=0$ where with $BA = \oplus_i (A[1])^{\otimes i}$
$$\WH{b}:BA \otimes M \otimes BA \to BA \otimes M \otimes BA$$
is an induced coalgebra map defined by
\begin{equation}
	\WH{b} \big( a_1 \otimes \cdots \otimes  a_i \otimes \UL{v} \otimes a_{i+1} \otimes \cdots \otimes a_{i+j} \big)
\end{equation}
$$= \sum (-1)^{K_1} a_1 \otimes \cdots \otimes m_{k}(a_l,\cdots,a_{l'}) \otimes \cdots \otimes a_i \otimes \UL{v} \otimes a_{i+1} \otimes \cdots \otimes a_{i+j} \big)$$
$$ +\sum (-1)^{K_2} a_1 \otimes \cdots \otimes b_{i-l+1,p-i}(a_l,\cdots,a_i, \UL{v}, a_{i+1} \otimes \cdots a_p) \otimes \cdots \otimes a_{i+j} \big)$$ 
$$+ \sum (-1)^{K_3}a_1 \otimes \cdots  a_i \otimes \UL{v} \otimes a_{i+1} \otimes \cdots \otimes m_{k}(a_q,\cdots,a_{q'}) \otimes \cdots\otimes a_{i+j} \big)$$
Here $K_*$ denotes the respective Koszul signs which occur when we move $m$ or $b$ through elements $a_*$ and $v$.
For example, it is easy to see that $K_1 = |a_1|' + \cdots + |a_{l-1}|'$.
We will underline an element of the module $M$ as above to distinguish it from elements of $A$.

For the case of $M=A[1]$, we may set $b_{k,l}= m_{k+l+1}$.
For the case of the dual $M = (A[1])^*[-n+2]$, we define $b^*_{k,l}$ as follows
\begin{equation}\label{def3}
b^*_{k,l}(x_1,\cdots,x_k,v^*,x_{k+1},\cdots,x_{k+l}) (w) =
(-1)^{\epsilon} v^* \big(  m_{k+l+1} (x_{k+1},\cdots,x_{k+l},w,x_1,\cdots,x_k) \big),
\end{equation}
where
$$\epsilon = 1+ K = 1+ |v^*|' + (|x_1|'+ \cdots + |x_k|')(|v^*|' + |x_{k+1}| + \cdots + |x_{k+1}|' + |w|').$$

One can check that this defines an $\AI$-bimodule structure, i.e. $\WH{b^*} \circ \WH{b^*} =0$.
The additional negative sign in $\epsilon$ is to cancel out the sign occurring from the switching of two $b$'s when we apply the
definition \ref{def3} twice in proving $\WH{b^*} \circ \WH{b^*} =0$. We leave the details to the reader.

The cyclic symmetric inner product can be understood as a special type of $\AI$-bimodule map $\psi:A \to A^*$.
Let $N$ be the degree of the inner product introduced in the last section.
Then, we are interested in the bimodule map $\psi:A \to A^*[-N]$, or
$\psi:A[1] \to (A[1])^*[-N+2]$  such that the bimodule map $\psi$  has degree zero.
Suppose $\{\psi_{k,l} \}$ defines an $\AI$-bimodule homomorphism $\psi:A \to A^*$.
\begin{lemma}\label{lemma:cyc}
Let $\psi$ be an $\AI$-bimodule homomorphism $\psi:A \to A^*$.
Define $$<a,b>=\psi_{0,0}(a)(b),$$
and suppose that $<,>$ is {\em non-degenerate}.
Then, it defines a cyclic symmetric inner product on $A$  if 
\begin{enumerate}\label{cycprop}
\item $ \psi_{k,l} \equiv 0$ for $(k,l) \neq (0,0)$
\item $\psi_{0,0}(a)(b) = -(-1)^{|a|'|b|'}\psi_{0,0}(b)(a).$
\end{enumerate}
Conversely, any cyclic symmetric inner product $<,>$ on $A$ give rise to an $\AI$-bimodule map
$\psi:A \to A^*$ with (1) and (2).
\end{lemma}
\begin{proof}
We first show `if' statement. Suppose the property (1) and (2). 
Then, $<,>$ is skew symmetric by (2). The cyclic symmetry of $<,>$ follows from 
the fact that  $\psi_{0,0}$ ( or $\psi$) defines an $\AI$-bimodule map along with (1) and (2) as follows.

An $\AI$-bimodule equation may be written in the cobimodule language as
$$\psi \circ \widehat{m} = m^* \circ \WH{\psi}.$$
Here we denote the $\AI$-bimodule structure of $A$ and $A^*$ as $m$ and $m^*$ respectively.
\begin{eqnarray*}
<m_k(a_1,\cdots,a_k),a_{k+1}> &=& \psi \big( m( \UL{a_1},a_2,\cdots,a_K) \big)(a_{K+1}) \\
&=& \psi \circ \WH{m}( \UL{a_1},a_2,\cdots,a_K)(a_{K+1}) \\
&=& m^* \circ \WH{\psi}( \UL{a_1},a_2,\cdots,a_K)(a_{K+1}) \\
&=& m^* (\UL{\psi_{0,0}(a_1)},a_2,\cdots,a_K)(a_{K+1}) \\
&=& -(-1)^{K_1} \psi_{0,0}(a_1) (m_{k}( a_2,\cdots,a_K,\UL{a_{K+1}})) \\
&=& +(-1)^{K_2} \psi_{0,0}(m_{k}( a_2,\cdots,a_K,\UL{a_{K+1}}))(a_1)  \\
&=& <m_k(a_2,\cdots,a_{k+1}),a_1>
\end{eqnarray*}
Here, $\psi$ is a map of degree zero, hence 
$$K_1 = |a_1|', K_2 = |a_1|' + |a_1|'(|a_2|' + \cdots +|a_{k+1}|').$$
The 1st, 3rd, 5th and 7th equality follows from the definitions, and the 2nd and 4th equality follows from (1).
This proves the first part of the lemma.
The second part of the lemma can be proved similarly.
%
\end{proof}

Before we proceed further, we recall some facts regarding $\AI$-bimodule structure.
Given an $\AI$-homomorphsm $f:A \to B$, it is known that $B$ has the structure of an $\AI$-bimodule over $A$. 
More precisely, let $m^A, m^B$ be the $\AI$-structures of $A$ and $B$ respectively.
And denote by $\WH{f}$ the induced coalgebra map from $f$.
Then, $\AI$-bimodule structure of $B$ over $A$ is given by
the maps
$b_{k,l}:A^{\otimes k} \otimes B \otimes A^{\otimes l} \to B$,
$$b_{k,l}(a_1 \otimes \cdots \otimes a_k \otimes b \otimes a_{k+1} \otimes \cdots \otimes a_{k+1})$$
$$= m^B \big(\WH{f}(a_1 \otimes \cdots \otimes a_k) \otimes b \otimes \WH{f}(a_{k+1} \otimes \cdots \otimes a_{k+1}) \big).$$
It is easy to check that $b_{k,l}$ defines an $\AI$-bimodule structure on $B$ over $A$.
Also, the $\AI$-homomorhism $f:A\to B$ now induces an $\AI$-bimodule homomorphism $\widetilde{f}:A\to B$ as a map between two $\AI$-modules over $A$.
In this case, $\widetilde{f}_{k,l}:A^{\otimes k} \otimes A \otimes A^{\otimes l} \to B $ is defined by $\widetilde{f}_{k,l} = f_{k+l+1}$.
One can check that $\widetilde{f}$ satisfies $ \widetilde{f} \circ \widehat{m}^A = m^B \circ \WH{\widetilde{f}}$.

Also, given an $\AI$-bimodule homomorphism $f:C \to D$, where $C$ and $D$ are $\AI$-bimodules over $A$, 
the $\AI$-bimodule map $f^*:D^* \to C^*$ is defined by 
\begin{equation}
	f^*_{k,l}(\vec{x},\UL{v},\vec{y})(w) = (-1)^{K} v \big(f_{l,k}(\vec{y},\UL{w},\vec{x}) \big)
\end{equation}
where $K = (\sum |x_*|')(|v|'+ \sum |y_*|'+|w|')$.

Now, we define the notion of a cyclic $\AI$-homomorphism as follows
\begin{definition}
Let $A,B$ be two cyclic symmetric $\AI$-algebras, hence with $\AI$-bimodule maps $\phi:A \to A^*$, $\psi:B \to B^*$
given by the previous lemma.
Then, an $\AI$-homomorphism $f:A \to B$ is called a {\em cyclic $\AI$-homomorphism } if
the following diagram of $\AI$-bimodules via $\AI$-bimodule homomorphisms commute.
\begin{equation}\label{comm1}
\xymatrix@C+1cm{ A \ar[r]^{\widetilde{f}}  \ar[d]_\phi^{cyc} & B \ar[d]_\psi^{cyc} \\ A^*  &  \ar[l]^{({\widetilde{f}})^*} B^*}	
\end{equation}
\end{definition}

The notion of a cyclic $\AI$-homomorphism was first defined by Kajiura (in the setting of non-commutative symplectic geometry
as maps preserving constant symplectic structures) and the above definition is in fact equivalent to it 
\begin{lemma}\label{cycmordef}$[$Definition 2.13 \cite{Kaj}$]$
An $\AI$-homomorphism $\{f_k\}_{k\geq 1 }$ between two $\AI$-algebras with cyclic symmetric inner products defines a cyclic $\AI$-homomorphism if and only if
\begin{enumerate}
\item $f_1$ preserves inner product $<a,b> = <f_1(a),f_1(b)>$.
\item \begin{equation}
\sum_{i+j=k, i,j>0} <f_i(x_1,\cdots,x_i), f_j(x_{i+1},\cdots,x_k)> =0.	
\end{equation}
\end{enumerate}
\end{lemma}
\begin{proof}
Proof can be easily seen from the commuting diagram above. More precisely, the first equation follows from commuting diagram of $\phi_{0,0}:A \to A^*$,
and the second equation follows from the vanishing of $\phi_{k,l}:A \to A^*$.
\end{proof}


Now, we define two homotopy notions of cyclic symmetric inner products.
The first one is given in the sense of $\AI$-bimodules, and the second one in the sense
of $\AI$-algebra, which can be easily seen by generalizing the commutative diagram above.

We fix an $\AI$-algebra $A$. First, we introduce the concept of a cyclic symmetric inner product of an $\AI$-bimodule over $A$.
\begin{definition}\label{bimodcyc}
The $\AI$-bimodule $M$ over $A$ is said to have a cyclic symmetric inner product if
for the $\AI$-bimodule differential $\{ b^M_{*,*} \}$,
we have non-degenerate skew symmetric bilinear map $<,> : M \times M \to \kk$ such that
\begin{equation}
	<b^M_{k,l}(a_1,\cdots,a_k,v,a_{k+1},\cdots,a_{k+l}),w> = (-1)^{K}<v,b^M_{l,k}(a_{k+1},\cdots,a_{k+l},w,a_1,\cdots,a_k)>
\end{equation}
or equivalently, we have a $\AI$-bimodule homomorphism $\psi:M \to M^*$ over $A$ with $\psi_{k,l} \equiv 0$ for $(k,l) \neq (0,0)$
with skew symmetric $\psi_{0,0}$.
Here $K = |v|' + (\sum_{i=1}^k |a_i|')(|v|' + \sum_{i=k+1}^{k+l} |a_i|' + |w|')$.
\end{definition}
Now, we define the notion of a homotopy inner product and the strong homotopy inner product.
\begin{definition}\label{hidef}
Let $C$ be an $\AI$-bimodule over $A$.
We call an $\AI$-bimodule map $\phi:C \to C^*$ a {\em homotopy inner product} if there exists an $\AI$-bimodule $D$ over $A$
which has a cyclic symmetric inner product $\psi : D \to D^*$ and a $\AI$-bimodule quasi-isomorphism $g:C \to D$ with
the commuting diagram
\begin{equation}\label{comm2}
\xymatrix@C+1cm{ C \ar[r]^g_{bimod}  \ar[d]_\phi & \,^\exists D \ar[d]_\psi^{cyc} \\ C^*  &  \ar[l]^{g^*} D^*}	
\end{equation}
\end{definition}
\begin{definition}\label{shidef}
Let $A$ be an $\AI$-algebra,
We call an $\AI$-bimodule map $\phi:A \to A^*$ a strong homotopy inner product if
there exists a cyclic symmetric $\AI$-algebra $B$ with $\psi:B \to B^*$ and
an $\AI$-quasi-isomorphism $f:A \to B$ such that the following diagram of $\AI$-bimodules over $A$ commutes 
\begin{equation}\label{comm3}
\xymatrix@C+1cm{ A \ar[r]_{ g = \widetilde{f}}  \ar[d]_\phi & \,^\exists B \ar[d]_\psi^{cyc} \\ A^*  &  \ar[l]^{g^*} B^*}	
\end{equation}
Here by $g:A \to B$, we denote the induced $\AI$-bimodule map $\widetilde{f}=g$ where $B$ is considered as an $\AI$-bimodule over $A$.
\end{definition}
Thus  $A$ is $\AI$-homotopy equivalent to $B$ which has a cyclic symmetric inner product, if and only if  $A$ has
a strong homotopy inner product by definition.

\section{Characterization of homotopy inner products}
Let $C$ be an $\AI$-bimodule over an $\AI$-algebra $A$.
For a given $\AI$-bimodule homomorphism $\phi:C \to C^*$, we find a condition on $\phi$ which is equivalent for $C$ to have a homotopy inner product.
(Definition \ref{hidef}).

\begin{theorem}\label{thm:hi}
An $\AI$-bimodule $C$ has a homotopy inner product if and only if 
there exists an $\AI$-bimodule map $\phi:C \to C^*$ satisfying the following skew-symmetry and homological non-degeneracy condition.
The skew symmetry condition is that 
\begin{equation}\label{bimodskew}
\phi_{k,l}(\vec{a},\UL{v}, \vec{b})(w) = -(-1)^{K}
\phi_{l,k}(\vec{b},\UL{w},\vec{a},v), 	
\end{equation}
for all $\vec{a} = (a_1,\cdots,a_k),\vec{b}=(b_1,\cdots,b_l)$ with $ a_i, v, b_j, w \in A$.
And we say $\phi$ is homologically non-degenerate, if for any non-zero $[a] \in H^*(C)$ with $a \in A$, there exists $[b] \in H^*(C)$ with $b \in C$
such that $\phi_{0,0}(a)(b) \neq 0$.
Here $K$ is the Koszul sign of switching $(\vec{a},v)$ with $(\vec{b},w)$:
$$K= (\sum |a_i|' + |v|')(\sum |b_i|' + |w|').$$
\end{theorem}
\begin{remark}
As we only require non-degeneracy on the homology level, it turns out that a homotopy inner product we find for $C$, say $\phi':C \to C^*$, is not
exactly the same as $\phi$, but only up to an $\AI$-bimodule quasi-morphism (in the 'if' statement). The map $\phi'$ will be constructed out of $\phi$.
If we have non-degeneracy on the chain level, we have $\phi = \phi'$.
\end{remark}
\begin{proof}
Let us first prove the 'only if' statement. By definition of a homotopy inner product,
$C$ is equipped with a map $\phi:C \to C^*$ which satisfy the commuting diagram \ref{comm2}.
Skew symmetry easily follows from the diagram as follows

Note that we have $\phi = g^* \circ  \WH{\psi} \circ \WH{g}$.
Hence, we have
$$\WH{g}(a_1\otimes \cdots \otimes a_k \otimes \UL{v} \otimes b_1 \otimes \cdots \otimes b_l) $$
$$= \sum_{i,j} a_1 \otimes \cdots \otimes a_{i} \otimes g_{k-i,j}(a_{i+1} \otimes \cdots \otimes a_k \otimes  \UL{v} \otimes b_1 \otimes \cdots \otimes b_j) 
\otimes b_{j+1} \otimes \cdots \otimes b_l) $$
and 
$$\phi_{k,l}(\vec{a},\UL{v}, \vec{b})(w) = \phi(a_1\otimes \cdots \otimes a_k \otimes \UL{v} \otimes b_1 \otimes \cdots \otimes b_l) (w)$$
$$=\sum_{i,j} g^* \big(a_1 \otimes \cdots \otimes a_{i} \otimes \psi_{0,0} (g_{k-i,j}(a_{i+1} \otimes \cdots \otimes a_k \otimes  \UL{v} \otimes b_1 \otimes \cdots \otimes b_j)) 
\otimes b_{j+1} \otimes \cdots \otimes b_l \big)(w)$$
$$=\sum_{i,j} (-1)^{K_1}<g_{k-i,j}(a_{i+1} \otimes \cdots \otimes a_k \otimes  \UL{v} \otimes b_1 \otimes \cdots \otimes b_j),
g_{l-j,i}(b_{j+1} \otimes \cdots \otimes b_l \otimes  \UL{w} \otimes a_1 \otimes \cdots \otimes a_{i})>, $$
$$=\sum_{i,j} -(-1)^{K_2}<
g_{l-j,i}(b_{j+1} \otimes \cdots \otimes b_l \otimes \UL{w} \otimes a_1 \otimes \cdots \otimes a_{i}),
g_{k-i,j}(a_{i+1} \otimes \cdots \otimes a_k \otimes  \UL{v} \otimes b_1 \otimes \cdots \otimes b_j)>, $$
$$= -(-1)^{K_3} \phi(b_1 \otimes \cdots \otimes b_l \otimes \UL{w} \otimes a_1\otimes \cdots \otimes a_k ) (v) = -(-1)^{K_3}
\phi_{l,k}(\vec{b},\UL{w},\vec{a},v),$$
where $K_3 = K$ in the skew symmetry condition and 
$$K_1 = (|a_1|'+ \cdots + |a_i|')(|a_{i+1}|' + \cdots + |a_k|'+|v|'+|b_1|'+\cdots+|b_l|'+|w|').$$
$$K_2 = K_1 + (|a_{j+1}|'+ \cdots + |v|' + |b_1|'+ \cdots + |b_j|')(|b_{j+1}|'+ \cdots + |b_l|' + |w|' + \cdots+|a_j|')$$
In the first equality, we used the fact that $\psi_{k,l} \equiv 0$ for $(k,l) \neq (0,0)$ and
in the third equality, we used the skew symmetry of $\psi$.
This proves the skew symmetry of $\phi$. 

For the homological non-degeneracy, first notice that inner product of $D$ is non-degenerate,
but it is not clear if it is homologically non-degenerate(with respect to $b_{0,0}$).
To show this, we first consider the Hodge type decomposition of $D = H \oplus S \oplus T$, where
$H$ is a subspace which contains each representative of the homology of $(D, b_{0,0})$ and $S = im (b_{0,0})$.
Cyclic symmetry implies that the inner product on $H \times S$ or $S \times H$ or $S \times S$ vanishes.
Hence, non-degeneracy implies that the inner product on $H \times H$ is non-degenerate, which 
implies the homological non-degeneracy.

Now, we prove the 'if' statement. We use the skew symmetric property to construct a new $\AI$-bimodule
$D$ ($\AI$-bimodule quasi-isomorphic to the original one) such that the $\AI$-bimodule $D$
has a cyclic symmetric inner product. Now, in general we will obtain the commuting diagram \ref{comm2},
not with $\phi$, but with a new map $\phi':C \to C^*$. But this still implies that $C$ has
a homotopy inner product, which finishes the proof.

We first consider the case that $\phi_{0,0}$ is non-degenerate on $C$ (
not only homologically, but also on the chain level). We will show
how to handle the general case at the end of the proof.

Here is how we proceed.
 As a vector space, we set $D=C$  and define $\psi :D \to D^*$ by setting
$\psi_{0,0} = \phi_{0,0}$ and set $\psi_{k,l} \equiv 0$ for $(k,l) \neq 0$. 
Instead of constructing $\AI$-bimodule structure on $D$, we first construct what is supposed to be an $\AI$-bimodule
map $g_{k,l} : C^{\otimes k} \otimes C \otimes C^{\otimes l} \to D$ satisfying the required equations (see the above or equation \ref{eq1} below) inductively. Secondly, we define an $\AI$-bimodule structure on $D$ inductively so that the already defined map
$g:C \to D$ becomes an $\AI$-bimodule homomorphism. The final step is to verify that 
the map $\psi:D \to D^*$ is indeed a map of $\AI$-bimodules, hence providing a cyclic symmetric inner product on $D$.
It will be clear from the construction that the diagram \ref{comm2} commutes. 

We begin our construction. Suppose that we are given an $\AI$-bimodule map $\phi:C \to C^*$ which is skew symmetric (\ref{bimodskew}). This implies that we have
$$\phi_{0,0}(v)(w) = (-1)^{|v|'|w|'+1} \phi_{0,0}(w)(v).$$ Hence so is $\psi_{0,0}$.
 
Now, we construct $g_{k,l}$ inductively. We set $g_{0,0}$ to be the identity map from $A$ to $A(=B)$.
Suppose we have constructed $g_{k,l}$ for all $k,l$ with $k+l < N$.
Consider any pair with $k+l =N$, and $a_i \in A$ for $1 \leq i \leq N$ and $v \in A$.

We determine $\phi_{k,l}(\vec{a},v,\vec{b})(w) \in \kk$ for $\vec{a} \in A^{\otimes k},\vec{b} \in A^{\otimes l},v \in C$. We recall that $\phi$ and $g$
are supposed to satisfy the following equations from the commuting diagram:
\begin{equation}\label{eq1}
\phi_{k,l}(\vec{a},v,\vec{b})(w) =
\end{equation}
$$\sum_{i,j} (-1)^{K_1} < g_{k-i,j}(a_{i+1},\cdots,a_{k},v,\cdots,a_{k+j}), g_{l-j,i}(a_{k+j+1},\cdots,a_{k+l},w,a_1,\cdots,a_i)>$$

Notice that by induction, the maps $g_{k-i,j}$ and $g_{i,l-j}$ have been already defined for $i>1$ or $j>1$.
The only undetermined expressions are the following two terms from the RHS of (\ref{eq1}).
\begin{equation}\label{eq2}
<g_{k,l}(\vec{a},v,\vec{b}),w> + (-1)^{K_2} < v,g_{l,k}(\vec{b},w,\vec{a})>
\end{equation}
Hence, considering all possible inputs $(\vec{a},v,\vec{b},w)$ among basis elements of $A$, we obtain linear system of equations, from
the non-degeneracy of $<,>$. (One can choose a good basis which is orthonormal or symplectic on suitable subspaces depending on the degree, so that the equation have a very simple form.)
 In fact, this linear system has a very simple structure as a linear system. Namely
the equation from data $(\vec{a},v,\vec{b},w)$ has the same unknown term as the equation from the data $(\vec{a}',v',\vec{b}',w')$
if and only if they are the same data or 
$$\vec{a} = \vec{b}',\vec{b} = \vec{a}', v = w',w=v'.$$
But the equation from the data $(\vec{a}',v',\vec{b}',w')$ in this case, exactly equals to the equation from the data $(\vec{a},v,\vec{b},w)$
up to the overall Koszul sign from the skew symmetry (\ref{bimodskew}). This is where we need skew symmetry condition.
(Without the skew symmetry assumption, the system would have no solution.)
Hence we can discard one of the equations from each pair, and solve each equation separately.

In fact, there exists a subtle problem when we solve the corresponding equation in the following particular case:
$$\vec{a}=\vec{b},\; v=w, \;\; k=l.$$
 Then, 
$$<g_{k,l}(\vec{a},v,\vec{b}),w> =<g_{l,k}(\vec{b},w,\vec{a}),v> = <g_{k,k}(\vec{a},v,\vec{a}),v>$$
In this case, the expression (\ref{eq2}) becomes
\begin{equation}\label{eq3}
	<g_{k,l}(\vec{a},v,\vec{a}),v> + (-1)^{(\sum |a_i|')(|v|'+\sum |a_i|'+|v|')}  <v, g_{k,l}(\vec{a},v,\vec{a})>
\end{equation}
The above expression could vanish in some case:
By comparing  this with the skew symmetry of the inner product
\begin{equation}
	<g_{k,l}(\vec{a},v,\vec{a}),v> + (-1)^{(\sum |a_i|'+ |v|'+\sum |a_i|')|v|'}  <v, g_{k,l}(\vec{a},v,\vec{a})>,
\end{equation}
one can see that the expression \ref{eq3} vanishes if  
$$ (\sum |a_i|')(|v|'+\sum |a_i|'+|v|') + (\sum |a_i|'+ |v|'+\sum |a_i|')|v|' $$
$$= \sum  |a_i|' + |v|' =0 (mod \, 2)$$
If this occurs, it may seem that there is no solution to the equation \ref{eq1}.
But, in this case, we show that all the other terms in the equation \ref{eq1} vanishes (as pairs) also, hence
the particular equation is satisfied trivially.
To see this, first consider the right hand side of the equation \ref{eq1}, which is $\phi_{k,k}(\vec{a},v,\vec{a})(v)$.
The vanishing of this term follows from the skew symmetry of $\phi$. Namely, 
\begin{eqnarray*}
\phi_{k,k}(\vec{a},v,\vec{a})(v) &=& -(-1)^{K}\phi_{k,k}((\vec{a},v,\vec{a})(v) \\
&=& - (-1)^{(\sum |a_i|' + |v|')(\sum |a_i|' + |v|')} \phi_{k,k}((\vec{a},v,\vec{a})(v) \\
&=& - \phi_{k,k}((\vec{a},v,\vec{a})(v)
\end{eqnarray*}
The last line follows from the assumption  $\sum  |a_i|' + |v|' =0$ above.

Also, other expressions in \ref{eq1} vanishes in pairs as follows.
we consider any two possible partitions of $\vec{a}$: $\vec{a} = (\vec{a^1},\vec{a^2})$ and  $\vec{a} = (\vec{a^3},\vec{a^4})$.
We denote by $|a^1|'$ the sum of degrees $\sum_{*} |a^1_{*}|'$.
Then, \ref{eq1} contains the following pair of expressions.
$$(-1)^{\epsilon_1}<g(\vec{a^2},v,\vec{a^3}),g(\vec{a^4},v,\vec{a^1})> 
+ (-1)^{\epsilon_2}<g(\vec{a^4},v,\vec{a^1}),g(\vec{a^2},v,\vec{a^3})> $$
where 
$$\epsilon_1 = |a^1|'(|a^2|'+|v|'+|a^3|'+|a^4|'+|v|'), \epsilon_1 = |a^3|'(|a^4|'+|v|'+|a^1|'+|a^2|'+|v|').$$
This term vanishes if the sign agrees with the following skew symmetry condition of $<,>$:
$$<g(\vec{a^2},v,\vec{a^3}),g(\vec{a^4},v,\vec{a^1})> + (-1)^{\epsilon_3}<g(\vec{a^4},v,\vec{a^1}),g(\vec{a^2},v,\vec{a^3})>=0 ,$$
where 
$$\epsilon_3= (|a^2|' + |v|' +|a^3|')(|a^4|' + |v|' +|a^1|')$$
Through elementary calculations, one can check that 
$$\epsilon_1+\epsilon_2 \equiv \epsilon_3 \,\, (mod \, 2),$$
under the assumption  $\sum  |a_i|' + |v|' =0$, which proves the claim.
In this way, we construct a map $g_{*,*}$ inductively which satisfies equations \ref{eq1}.
And by construction, $\phi = g^* \circ \psi \circ g$.
To finish the construction of the new $\AI$-bimodule, we construct a new bimodule structure $b^D_{*,*}$ on $D(=C)$ such that
$g:C \to D$ and $\psi(=\psi_{0,0}):D \to D^*$  are maps between $\AI$-bimodules.
First, we choose an $\AI$-bimodule structure so that $g$ becomes a bimodule map.
\begin{lemma}
Let $(C, b^C_{*,*})$ be an $\AI$-bimodule over an $\AI$-algebra $A$  and consider
a family of bilinear maps $$g_{i,j} : A^{\otimes i} \otimes C \otimes  A^{\otimes i} \to D.$$ for
$i,j \in \NN \cup \{0\}$, and $g_{0,0}:C \to D$ is a vector space isomorphism.
Then, there exists a new  bimodule structure $b^D_{*,*}$ on $D(=C)$  such that $g$ becomes
an $\AI$-bimodule homomorphism between these two $\AI$-bimodules.
\end{lemma}
\begin{proof}
The new bimodule structure is defined inductively 
using the $\AI$-bimodule equation as follows.
We would like to have
\begin{equation}\label{eq5}
 (g \circ \WH{b}^C)(\vec{x} \otimes v \otimes \vec{y}) = (b^D \circ \WH{g}) (\vec{x} \otimes v \otimes \vec{y}).
\end{equation}
Suppose $b^D_{k,l}$ is defined for $k+l <N$.
Consider the case $k+l=N$.
Then, in the equation \ref{eq5}, the terms on the left hand side is already given, and the terms on the right hand side
is determined by induction hypothesis except the term $b_{k,l}(\vec{x},g_{0,0}(v),\vec{y})$.
Hence, this term is uniquely determined by the equation \ref{eq5}.

To show that the maps $b_{k,l}^D$ indeed gives an $\AI$-bimodule structure on $D$, we need to show that
$\WH{b}^D \circ \WH{b}^D =0.$
But, first, it is not hard to note that the map
$$\WH{g} : BA \otimes C \otimes BA \to BA \otimes D \otimes BA.$$
is a surjective since $g_{0,0}$ is a vector space isomorphism. 

Hence it is enough to check $ \WH{b}^D \circ \WH{b}^D \circ \WH{g} =0.$
But since $\WH{b}^D \circ \WH{g} = \WH{g} \circ \WH{b}^C$,
$$ \WH{b}^D \circ \WH{b}^D \circ \WH{g} = \WH{g} \circ \WH{b}^C \circ \WH{b}^C =0.$$
\end{proof}

Now, we show that $\psi:D \to D^*$ is an  $\AI$-bimodule homomorphism with the bimodule structure $b^D$ on $D$ and  $b^{D^*}$ on $D^*$. Namely, we need to prove that 
\begin{equation}\label{eq6}
\WH{\psi} \circ \WH{b}^D = \WH{b}^{D^*} \circ \WH{\psi}.
\end{equation}

The map $\phi$ is an $\AI$-bimodule map, hence 
$$\WH{\phi} \circ \WH{b}^C = \WH{b}^{C^*} \circ \WH{\phi}.$$
Hence, using $\WH{\phi} =\WH{g}^* \circ \WH{\psi} \circ \WH{g}$,
$$\WH{g}^* \circ \WH{\psi} \circ \WH{g} \circ \WH{b}^C = \WH{b}^{C^*} \circ \WH{g}^* \circ \WH{\psi} \circ \WH{g}.$$
Therefore, we have
$$\WH{g}^* \circ \big( \WH{\psi}  \circ \WH{b}^D \big) \circ \WH{g} = \WH{g}^* \circ \big( \WH{b}^{D^*} \circ \WH{\psi} \big) \circ \WH{g}.$$
This in fact proves the equality \ref{eq6}, because $\WH{g}$ is surjective.
This finishes the construction of $D$, hence gives the proof of the theorem assuming the non-degeneracy on the chain level.

In the general case where $C$ has only homological non-degenerate inner product $\phi_{0,0}$,
we proceed as follows.
By the minimal model theorem(see \cite{FOOO} or \cite{TZ}), there exists an $\AI$-bimodule $H^*(C)$ over $A$, with
an $\AI$-bimodule quasi-morphism $\eta : H^*(C) \to C$. Thus, we obtain the map
$\widetilde{\phi}$ from the
following commuting diagram.
\begin{equation}\label{comm4}
\xymatrix@C+1cm{ C   \ar[d]_\phi  &  \ar[l]^\eta  H^*(C) \ar[d]_{\widetilde{\phi}} \\ C^*  \ar[r]^{\eta^*} &   (H^*(C))^*}	
\end{equation}
Note that the induced map $\widetilde{\phi}_{0,0}: H^*(C) \to (H^*(C))^* $ is non-degenerate.
Also, it is not hard to show that if $\phi$ satisfies skew symmetry, then $\widetilde{\phi}$ also satisfies the
skew symmetry.
Now, we can apply the construction above in the non-degenerate case for the map $\widetilde{\phi}_{0,0}: H^*(C) \to (H^*(C))^* $.
Therefore, there exists a bimodule $D$, which has cyclic symmetric inner product $\psi:D \to D^*$, and an $\AI$-bimodule map $g:H^*(C) \to D$ with the commuting diamgram \ref{comm5}.

But since $\eta$ is only a quasi-isomorphism, we cannot invert the arrow of $\eta$ exactly.
Namely, there exists an quasi-inverse $\xi$ of $\eta$ so that $\eta \circ \xi$ is only homotopic to identity.
So, as an alternative,  by using $\xi$, we construct a new map $\phi':C \to C^*$ from $\widetilde{\phi}$ using the following diagram
\begin{equation}\label{comm5}
\xymatrix@C+1cm{ C  \ar[r]^\xi \ar@{.>}[d]_{\phi'}  &  \ar[r]^g  H^*(C) \ar[d]_{\widetilde{\phi}} & D \ar[d]_{\psi} 
 \\ C^*  &   \ar[l]^{\xi^*} (H^*(C))^*	& \ar[l]^{g^*} D^* }
\end{equation}
(It is easy to see that $\phi$ and $\phi'$ are related by $\eta \circ \xi$). 
The above diagram proves that $C$ has a homotopy inner product $\phi':C \to C^*$.
\end{proof}

\section{Characterization of strong homotopy inner products}
Let $A$ be an $\AI$-algebra, and consider the induced $\AI$-bimodule structures on $A$ and $A^*$.
For a given $\AI$- bimodule homomorphism $\phi:A \to A^*$, we find a condition on $\phi$ which is equivalent for $A$ to have
a strong homotopy inner product.

\begin{theorem}\label{thm:sh}
An $\AI$-algebra $A$ has a strong homotopy inner product if and only if 
there exists an $\AI$-bimodule map $\phi:A \to A^*$, satisfying the following three conditions.
\begin{enumerate}
\item (Skew symmetry)
$$\phi_{k,l}(\vec{a},v, \vec{b})(w) = -(-1)^{K}\phi_{l,k}(\vec{b},w,\vec{a},v).$$
\item (Closedness)
For any choice of a family $(a_1,\cdots,a_{l+1})$ and any choice of indices $1 \leq i<j<k\leq l+1$,
we have
$$(-1)^{K_i} \phi(.. ,\UL{a_i},..)(a_j) + (-1)^{K_j} \phi(.. ,\UL{a_j},..)(a_k) +
(-1)^{K_k} \phi(.. ,\UL{a_k},..)(a_i) =0,$$
where the arguments inside $\phi$ are uniquely given by the cyclic order of the family $(a_1,\cdots,a_{l+1})$,
and the signs $K_*$ are given by the Koszul convention:
\begin{equation}\label{eq}
K_* = (|a_1|' + \cdots + |a_*|')(|a_{*+1}|'+\cdots + |a_k|').
\end{equation}
\item (Homological non-degeneracy) For any non-zero $[a] \in H^*(A)$ with $a \in A$, there exists a $[b] \in H^*(A)$ with $b \in A$,  such that
$\phi_{0,0}(a)(b) \neq 0$.
\end{enumerate}
\end{theorem}
\begin{remark}
As in the case of $\AI$-bimodule, a newly constructed map $\phi':A \to A^*$ ( not the map $\phi$), becomes a strong homotopy inner product
with the homological non-degeneracy condition. If it is non-degenerate on $A$, then $\phi'=\phi$.
\end{remark}
\begin{proof}
We first show that a strong homotopy inner product, $\phi$, satisfies these properties.
The first property, Skew symmetry, follows from the similar arguments in the previous section due to the relation $b_{k,l}=m_{k+l+1}$.
Also note that $\AI$-homomorphisms $f:A \to B$ induces an $\AI$-bimodule
 map $g:A \to B$ and they are related by $g_{k,l}=f_{k+l+1}$.

To show the closedness condition, we apply the equation \ref{eq1} to each term of the equation. More precisely, the first term is 
\begin{eqnarray*}
(-1)^{K_i}\phi(.. ,\UL{a_i},..)(a_j)&=& (-1)^{K_i}\phi(a_{j+1},a_{j+2},\cdot, \UL{a_i},\cdots,a_{j-1})(a_j) \\
&=& \sum (-1)^{K_1}<f(\cdots,a_i,\cdots),f(\cdots,a_j,\cdots)> \\
&=& \sum (-1)^{K_2}<f(\cdots,a_k,\cdots ,a_i,\cdots),f(\cdots,a_j,\cdots)>\\
 &&+ \sum (-1)^{K_3}<f(\cdots,a_i,\cdots),f(\cdots,a_j,\cdots,a_k,\cdots)>
\end{eqnarray*}
Similarly, 
\begin{eqnarray*}
(-1)^{K_j}\phi(.. ,\UL{a_j},..)(a_k) &=& \sum (-1)^{K_4}<f(\cdots,a_i,\cdots ,a_j,\cdots),f(\cdots,a_k,\cdots)>\\
&&+ \sum (-1)^{K_5}<f(\cdots,a_j,\cdots),f(\cdots,a_k,\cdots,a_i,\cdots)>,\\
(-1)^{K_k}\phi(.. ,\UL{a_k},..)(a_i) &=& \sum (-1)^{K_6}<f(\cdots,a_j,\cdots ,a_k,\cdots),f(\cdots,a_i,\cdots)> \\
&&+ (-1)^{K_7}\sum <f(\cdots,a_k,\cdots),f(\cdots,a_i,\cdots,a_j,\cdots)>,
\end{eqnarray*}
It is easy to note that each expression in the overall sum, occurs in pairs such as
$$(-1)^{K_8}<f(\cdots,a_k,\cdots ,a_i,\cdots),f(\cdots,a_j,\cdots)>+ $$
$$(-1)^{K_9}<f(\cdots,a_j,\cdots),f(\cdots,a_k,\cdots,a_i,\cdots)>,$$
which cancel out due to the skew symmetry condition as all the related signs arise from the Koszul convention. This proves that the closedness condition
is satisfied for a strong homotopy inner product. 
Homological non-degeneracy can be proved similarly as in the case of homotopy inner products.

Now, we prove the 'if' statement. Suppose $\phi:A \to A^*$ satisfies the three conditions of the theorem.
First, we assume that $\phi_{0,0}$ is non-degenerate on $A$, not only homologically non-degenerate 
as in the previous section. General case will be discussed at the end of the proof.
We construct an $\AI$-algebra $B$ with cyclic symmetry, which is $\AI$-homotopy equivalent to $A$.

This case is more complicated than the $\AI$-bimodule case, due to the abundance of the equations to be satisfied.
Our strategy is to use only a small part of the bimodule map , $\phi_{0,*}$ to construct an $\AI$-algebra $B$ and a map $f:A \to B$ and show that
the compatibility with $\phi_{*,*}$ automatically follow from the closedness condition.

We set $B=A$ as a vector space, and  we will construct maps $f_k:A^{\otimes k} \to B$ inductively. We set $f_1$ to be the identity map, and
$\psi_{0,0} =\phi_{0,0}$.

The commuting diagram \ref{comm3} gives rise to equations (cf. equation \ref{eq1}) for $a_* \in A$.
\begin{equation}\label{seq1}
	\phi_{0,k-1}(\UL{a_1},a_2,\cdots,a_k)(a_{k+1})
\end{equation}
$$= \sum_i <f_i(a_1,\cdots,a_i),f_{k-i+1}(a_{i+1},\cdots,a_{k+1})>$$

Now, we will construct family of maps $f_{k}$ inductively satisfying the above equation \ref{seq1}.
Suppose we have determined $f_i$ for $i < k$.
In the equation \ref{seq1}, all the terms are already determined except the following two terms from the right hand side
$$<f_k(a_1,\cdots,a_k),a_{k+1}> + < a_1,f(a_2,\cdots,a_{k+1})>$$
By considering a suitable good basis of $A$ with respective to the inner product, the equations \ref{seq1} provide a system of linear equations as $a_i$ varies over basis elements. This linear system in fact can be partitioned into several linear subsystems.
This is because the expression $<f_k(a_1,\cdots,a_k),a_{k+1}>$  appears exactly twice in the linear system:
Namely, consider the equation \ref{seq1} with $\phi_{0,k-1}(\UL{a_{k+1}},a_1,\cdots,a_{k-1})(a_{k})$. Then
the equation contains the expression $<a_{k+1},f_k(a_1,\cdots,a_k)>$ which equals $(-1)^{K}<f_k(a_1,\cdots,a_k),a_{k+1}>$.
In fact, a moments thought tells us that we need to consider a linear subsystem which is obtained by the equations as we cyclically rotate the arguments $(a_1,\cdots,a_{k+1})$ 
 
Thus obtained linear subsystem is now somewhat similar to the following simple example:
$$x_1 - x_2 = c_1, x_2-x_3 =c_2, \cdots,  x_{k+1}-x_1 = c_{k+1}.$$
This linear system has a solution only if $ \sum_j c_j =0$, since the coefficient system is not non-degenerate.

Back to our original problem, we also need to check certain identities to make sure the existence of a solution.
First, we consider the case that all  elements $a_1, \cdots,a_{k+1}$ are distinct from each other.
Then, consider the sum over all cyclic permutations: For $\sigma \in \ZZ/ (k+1)\ZZ$, (with $a_i=a_j$ if $i=j$ mod $k+1$)
\begin{equation}\label{seq2}
	\sum_{\sigma} (-1)^{K_\sigma}\phi_{0,k-1}(\UL{a_{\sigma+1}},a_{\sigma+ 2},\cdots,a_{\sigma+k})(a_{\sigma+k+1})
\end{equation}
\begin{equation}\label{seq3}
= \sum_{\sigma} (-1)^{K_\sigma} \sum_i <f_i(a_{\sigma+1},\cdots,a_{\sigma+i}),f_{k-i+1}(a_{{\sigma+i+1}},\cdots,a_{{\sigma+k+1}})>	
\end{equation}
where $$K_\sigma = (|a_1|' + \cdots + |a_{\sigma-1}|')(|a_{\sigma}|' + \cdots + |a_{k+1}|')$$
In the above expression \ref{seq3}, each term appears in pairs from a different permutation
$$(-1)^{K_1} <f_i(a_{\sigma+1},\cdots,a_{\sigma+i}),f_{k-i+1}(a_{{\sigma+i+1}},\cdots,a_{{\sigma+k+1}})>	$$
$$+ (-1)^{K_2} <f_{k-i+1}(a_{{\sigma+i+1}},\cdots,a_{{\sigma+k+1}}),f_i(a_{\sigma+1},\cdots,a_{\sigma+i}),>$$
and these terms vanish due to the skew symmetry of $<,>$.
This implies that the equation (\ref{seq2},\ref{seq3}) maybe considered as a redundant equation of the system, and this system has
a solution if the expression \ref{seq2} vanishes. We will show this vanishing from the closed condition as follows.

\begin{lemma} The expression \ref{seq2} vanishes.
\end{lemma}
\begin{proof}
We fix the family  $(a_1,\cdots, a_{k+1})$ and do not change throughout the proof of the lemma.
We may write (for $i \neq j$) 
$$[a_i,a_j] = (-1)^{K_i} \phi_{i+k-j,j-i}(\cdots,\UL{a_i},\cdots)(a_{j})$$
and the skew symmetry implies that  
 $$[a_i,a_j]=-[a_j,a_i].$$
Also, the closed condition reads as
$$[a_i,a_j] + [a_j,a_k] + [a_k,a_i]=0.$$
And, the statement we need to prove can be written as
\begin{equation}\label{seq4}
[a_1,a_{k+1}] + [a_2,a_1] + [a_3,a_2] + \cdots + [a_{k+1},a_{k}]=0
\end{equation}
But this can be proved easily from the above two equations.
%
\end{proof}

This proves that each linear subsystem has a solution when all $a_i$'s are distinct. 

The case when all $a_i$'s are not distinct can be dealt analogously.
The case we need to be careful is when for some $\sigma \in  \ZZ/ (k+1)\ZZ,\;(\sigma_0 \neq 0)$,
$$(a_1,\cdots,a_{k+1}) = (a_{\sigma+1},a_{\sigma+2},\cdots,a_{\sigma+k+1}).$$
In the case that  $\sigma \neq 1$, 
 instead of considering all cyclic permutation in \ref{seq2}, we consider the 
cyclic permutations corresponding to $0,1,\cdots,\sigma-1 \in \ZZ/ (k+1)\ZZ$ and we can proceed similarly.

The case $\sigma=1$ is more subtle. In this case, 
we have a family $\vec{a} =(a,a,\cdots,a)$ for $a \in A$, and a single equation:
\begin{equation}\label{aaa}
\phi_{0,k-1}(\UL{a},a,\cdots,a)(a)=
\sum<f_i(a,\cdots,a),f_{k-i+1}(a,\cdots,a)>.
\end{equation}
In this case, like in the case of $\AI$-bimodule, two undetermined terms in the inductive process, may
cancel out each other. More precisely, from the skew symmetry
$$<a,f(a,\cdots,a)> + <f(a,\cdots,a),a> = (1 - (-1)^{|a|'(k|a|')} <f(a,\cdots,a),a>,$$
the undetermined two terms (LHS) cancel each other out if $k|a|'$ is even.
We show that in this case, the whole equation \ref{aaa} actually vanishes.

First, if $k|a|'$ is even, we show that the RHS of the equation \ref{aaa} vanishes.
Terms in the RHS comes in pairs
$$<f_i(a,\cdots,a),f_{k-i+1}(a,\cdots,a)> + <f_{k-i+1}(a,\cdots,a),f_i(a,\cdots,a)>.$$
And from the skew symmetry condition, this terms cancels because
$$ k|a|' \;\textrm{even} \; \Rightarrow \;(-1)^{(i|a|')(k-i+1)|a|'} =1. $$
If $|a|'$ is even, this is clear. If $k$ is even, $i$ and $k-i+1$ have different parity.

Thus, we only need to show that the LHS of the equation \ref{aaa},
$\phi_{0,k-1}(\UL{a},a,\cdots,a)(a)$ vanishes. To show this, we use
the skew symmetry and the closed condition of $\phi$. In fact,
we will use the equation \ref{seq4} which was derived from these two
properties. Namely, we have
\begin{equation}\label{aaa2} 
[1,k+1] + [k+1,k] + \cdots + [3,2] + [2,1]=0.
\end{equation}
In this case, since all $a_i$'s are equal, 
Hence, $$[i+1,i] = (-1)^{k|a|'} [i,i-1].$$
Therefore, the equation \ref{aaa2} is equivalent to
$(k+1)[2,1] =0.$
Hence, $[2,1]$ vanishes, and so does $[1,2]$.
This implies the claim. 
This finishes the construction of $f$.

Now, one can construct an $\AI$-algebra structure on $B$ which is homotopy equivalent to $A$ such that
$f:A \to B$ is an $\AI$-algebra homomorphism.
\begin{lemma}
Let $(A, m)$ be an $\AI$-algebra and consider
a family of bilinear maps $f_{k} : A^{\otimes k} \to B.$ for
$k \in \NN $, and $f_1$ is a vector space isomorphism.
Then, there exists a $\AI$-algebra structure on $B(=A)$  such that $f$ defines
an $\AI$-algebra quasi-isomorphism between these two $\AI$-algebras.
\end{lemma}
\begin{proof}
The proof is analogous to the case of $\AI$-bimodule and we omit its proof.
\end{proof}

The $\AI$-homomorphism $f:A \to B$ induces an $\AI$-bimodule map $g:A \to B$ where $B$ is regarded as an $\AI$-bimodule over $A$.
Then, from the composition of the following 
$$A \stackrel{g}{\rightarrow} B \stackrel{\psi}{\rightarrow} B^* \stackrel{g^*}{\rightarrow} A^*,$$
we have a map $\widetilde{\phi}_{k,l}:A \to A^*$ $(\widetilde{\phi} = g^* \circ \psi \circ g)$.
By construction, we have $$\widetilde{\phi}_{0,k} = \phi_{0,k}.$$
Note that we do not know whether $\psi:B \to B^*$ is a map of $\AI$-bimodules yet, and the composed map 
$\widetilde{\phi}_{k,l}:A \to A^*$ is not an $\AI$-bimodule map yet. Whether it is a bimodule map depends on the relation of $\widetilde{\phi}$ and $m$,
and we only consider the map $\widetilde{\phi}$ itself.

Now, we show that $\widetilde{\phi}$ actually is the same as $\phi$.For this, note that from the diagram and the arguments in the beginning of proof
of the theorem,  $\widetilde{\phi}$ also satisfies the skew symmetry and the closed condition. 
Using this, we prove that 
\begin{lemma} For any $k,l \in \NN \cup \{0 \}$, we have 
$$\widetilde{\phi}_{k,l} = \phi_{k,l}$$
\end{lemma}
\begin{proof}
It is enough to prove that $\phi_{k,l}$ (or $\widetilde{\phi}_{k,l}$) can be written as a linear combination of $\phi_{0,*}$'s.
We fix the family of elements $(a_1,\cdots, a_{k+l+1})$ and do not change throughout the proof of the lemma.
Note that, we already have 
$$\phi_{0,k+l+1}(\UL{a_1},\cdots,a_{k+l})( a_{k+l+1}) = \widetilde{\phi}_{0,k+l+1}(\UL{a_1},\cdots,a_{k+l})( a_{k+l+1}).$$ 
Notice that  the term $\phi_{0,k+l+1}(\UL{a_i},\cdots,a_{i-2})( a_{i-1})$ may be written as $[a_i,a_{i-1}]$.
And 
$$\phi_{k,l}(a_{j+1},\cdots,\UL{a_i},a_{i+1},\cdots,a_{j-1})(a_j)= [a_i,a_j].$$
Then, the lemma can be easily proved by the identity
$$[a_i,a_j] = [a_i,a_{i+1}] + [a_{i+1},a_{i+2}] + \cdots + [a_{j-1},a_j]$$
\end{proof}

Thus we have a commuting diagram \ref{comm3}.
Thus it remains to show that $\psi:B \to B^*$ is a map of $\AI$-bimodules, which can be proved exactly as in the case of
$\AI$-bimodules using the commuting diagram.
This finishes the proof of the theorem in the case that the $\phi_{0,0}$ is non-degenerate.

For the case that $\phi_{0,0}$ is homologically non-degenerate, we proceed as follows.
The following lemma is easy to check.
\begin{lemma}
Given an $\AI$-homomorphism $h:P \to Q$, and an $\AI$-bimodule map $\chi:Q \to Q^*$ between two $\AI$-bimodules over $Q$, 
there exists an $\AI$-bimodule map $\zeta:P \to P^*$ with the commuting diagram of $\AI$-bimodules over $P$:
\begin{equation}\label{comm7}
\xymatrix@C+1cm{ P \ar[r]_{ \widetilde{h}}  \ar@{.>}[d]_\zeta &  Q  \ar[d]_{\widetilde{\chi}} \\ P^*  &  \ar[l]^{\widetilde{h}^*} Q^*}	
\end{equation}
Here, $\widetilde{\chi}:Q \to Q^*$ is an $\AI$-bimodule map over an $\AI$-algebra $P$, which is induced from $\chi$.

Moreover, if $\chi$ satisfies the skew symmetry and closedness condition, so is $\zeta$.
\end{lemma}

By the minimal model theorem for an $\AI$-algebra $A$, there exists an induced $\AI$-structure on $H^*(A)$ such that
there exists an $\AI$-quasi-isomorphism $h : H^*(A) \to A$. Therefore, by the lemma, 
we obtain an $\AI$-bimodule map $\zeta: H^*(A) \to (H^*(A))^*$ which satisfies skew symmetry, closedness condition, and
non-degeneracy. Now, we can use the construction so far to define an $\AI$-algebra $B$ with the cyclic symmetric
inner product $\psi:B \to B^*$, and an $\AI$-homomorphism $f:H^*(A) \to B$ with the commuting diagram \ref{comm8}
of $\AI$-bimodules over $H^*(A)$.
\begin{equation}\label{comm8}
\xymatrix@C+1cm{  \ar[r]^{g=\widetilde{f}}  H^*(A) \ar[d]_{\zeta} & B \ar[d]_{\psi} 
 \\  (H^*(A))^*	& \ar[l]^{g^*} B^* }
\end{equation}

But since $h$ is only a quasi-isomorphism, we cannot invert the arrow of $h$ exactly.
There exists only a quasi inverse $w:A \to H^*(A)$ with $h \circ w$ is only homotopic to the identity map.
So, as an alternative,  by using $w$, we construct a new map $\phi':A \to A^*$ from $\zeta$ using the following diagram
and the above lemma. (Note that the above diagram \ref{comm8} induces the same diagram
of $\AI$-bimodules over $A$, using the quasi-isomorphism $w$, where we denote the induced maps using $\widetilde{\,}\;$)
\begin{equation}\label{comm9}
\xymatrix@C+1cm{ A  \ar[r]^w \ar@{.>}[d]_{\phi'}  &  \ar[r]^{\widetilde{g}}  H^*(A) \ar[d]_{\widetilde{\zeta}} & B \ar[d]_{\widetilde{\psi}} 
 \\ A^*  &   \ar[l]^{w^*} (H^*(A))^*	& \ar[l]^{\widetilde{g^*}} B^* }
\end{equation}
This proves that $A$ has a strong homotopy inner product $\phi'$.
\end{proof}

\section{Relation to the non-constant symplectic form}\label{sec:nonc}
According to Kontsevich (see \cite{KS},\cite{K},\cite{Kaj}), $\AI$-algebra, say $A$, may be considered as a formal non-commutative supermanifold, say $X$. If $A$ has basis elements $\{e_i\}$, then, consider the formal dual coordinates $x_i$. 
In this formalism, operations $\{m_*\}$ correspond to the vector field $\delta$ on $X$.
If $$m_k(e_{i_1},\cdots,e_{i_k}) = \sum_j m^{j}_{i_1,\cdots,i_k} e_{j}$$
Then, we have (with Einstein convention)
$$ \delta = \sum_{k=1}^\infty \overleftarrow{\frac{\partial}{\partial x_j} } m^{j}_{i_1,\cdots,i_k} x_{i_k}
\cdots x_{i_1}.$$
A cyclic symmetric inner product corresponds to a (constant) symplectic structure $\omega$ on $X$.
Namely, let $<e_i,e_j> = \omega_{i,j}$ and consider its inverse $\omega^{i,j}$
Then, $$(\;\;,\;\;):= \overleftarrow{\frac{\partial}{\partial x_i}} \omega^{i,j} 
\overrightarrow{\frac{\partial}{\partial x_j}}$$
defines a Poisson bracket (See \cite{Kaj} for detailed explanations).

As in the section \ref{sec:po}, cyclic $\AI$-algebra is equipped with a cyclic potential $\Phi$.
There exist a relation $\delta = (\, , \Phi)$, and on the other hand, any vector field ($\AI$-structure)
$\delta$ obtained this way has a cyclic symmetry.

Kajiura has extended  the notion of the constant symplectic structure to that of a non-constant symplectic structures,
and showed that any non-constant symplectic form on a formal non-commutative supermanifold can be transformed to
be constant by a coordinate transform. The definition of non-constant symplectic form actually came from
a physical motivation via open string theory, and it was defined in Definition 4.5 of \cite{Kaj} as a poisson bracket
$$(A,B) = \sum_{i,j,I,J} (-1)^{(B-j)J}\omega^{ij}_{JI} \big((A \overleftarrow{\frac{\partial}{\partial x_i}})
x^I (\overrightarrow{\frac{\partial}{\partial x_j}}B) x^J \big)_c$$
which satisfies the skew symmetry and the Jacobi identity,
or as a non-constant symplectic two form $$\omega = \sum_{i,j,I,J} \omega_{ji,JI}(x^I dx^i x^J dx^j)_c$$
where the coefficient satisfies 
$$\omega_{ji,JI} = -\omega_{ij,IJ},$$ and the two form being closed.

\begin{theorem}\label{thm:equiv}
Let $A$ be an $\AI$-algebra which is finite dimensional.
The strong homotopy inner product on $A$ with non-degeneracy condition (not only homologically non-degenerate) is equivalent to a non-constant symplectic form on the corresponding formal non-commutative supermanifold $X$.
\end{theorem}
\begin{proof}
Intuitively, this is clear. A strong homotopy inner product, can be made into a strict cyclic symmetric inner 
product after a coordinate change (via $\AI$-quasi-isomorphism). Also, a non-constant symplectic structure
can be made into the constant symplectic structure via non-commutative version of the Darboux theorem
(\cite{Kaj},\cite{G}).

The characterization  theorem of a  strong homotopy inner product provides a rigorous proof of this correspondence.
Namely, for the basis $\{e_j\}$ of $A$, let $\phi:A \to A^*$ be a strong homotopy inner product.
Then, for $I = \{i_1,\cdots,i_a \}, J=\{j_1,\cdots,j_b \}$,
 we set 
 $$\omega_{ij,IJ}= \phi(e_{i_1},\cdots,e_{i_a},\UL{e_i},e_{j_1},\cdots,e_{j_b})(e_j)$$
Then, the skew symmetry of $\omega$ , and that of $\phi$  in the Theorem \ref{thm:sh} exactly corresponds to each other, and
the closedness condition of $\omega$ and that of $\phi$ in the Theorem \ref{thm:sh} also exactly corresponds to each other. Hence the proposition follows from the Theorem \ref{thm:sh}.
\end{proof}

\section{Potential of a cyclic $\AI$-algebra and its invariance}\label{sec:po}
In this section, we introduce Kajiura's work on the invariance of the symplectic action under the
cyclic homomorphism, which is not well-known to the mathematical community. If we translate it into a mathematical language, it means that cyclic $\AI$-algebra $A$ is equipped
with a potential whose main feature is the invariance under cyclic $\AI$-homotopy equivalence up to reparametrization and cyclization.
Let $(A,m^A_*)$ be a cyclic $\AI$-algebra.
Let $e_i$ be generators of $A$ as a vector space, which is assumed to be finite dimensional.
We may transfer the $\AI$-structure on $A$ to $H^*(A)$ to get a finite dimensional $\AI$-algebra.
Define $\bx = \sum_i e_i x_i$ where $x_i$ are formal parameters with $deg(x_i) = - deg(e_i)$.

\begin{definition}\label{podef}
We define 
\begin{equation}
\Phi^A(\bx) = \sum_{k=1}^\infty  \frac{1}{k+1} < m^A_k( \bx,\bx,\cdots,\bx),\bx>	
\end{equation}
\end{definition}
This may be considered as a systematic way of gathering structure constants of a cyclic $\AI$-algebra.
Although two homotopy equivalent $\AI$-algebras are interwined by an $\AI$-homomorphism, 
this is not good enough to interwine their structure constants. One of the benefit of having
cyclic $\AI$-algebra and cyclic morphism is  this invariance of structure constants up
to the reparametrization.

We first set up some notations. We write (with Einstein convention)
\begin{equation}
m^A_k(e_{i_1},\cdots,e_{i_k}) = \,^Am^i_{i_1,\cdots,i_k} e_{i} ,\;\;\; h^i_{j_1,\cdots,j_k} \in \kk
\end{equation}

Consider a cyclic $\AI$-algebra $(B,m^B_*)$ and 
a cyclic $\AI$-homomorphism $h:B \to A$. We define the change of coordinate as follows.
We also assume $B$ is finite dimensional as a vector space, and denote by $\{f_*\}$ its basis, and 
introduce corresponding formal variables $y_*$ as before.
Suppose 
$$h_{k}(f_{j_1},\cdots,f_{j_k}) = h^i_{j_1,\cdots,j_k} e_i, \;\;\; h^i_{j_1,\cdots,j_k} \in R,$$
where $R$ is the ring for the $\AI$-algebras.
Then, we set
\begin{equation}\label{changeco}
	x_i \mapsto h^i_{j_{11}} y_{j_{11}} + h^i_{j_{21},j_{22}} y_{j_{21}} y_{j_{22}} + \cdots +
	h^{i}_{j_{l1},\cdots,j_{lk}} y_{j_{l1}}\cdots y_{j_{lk}} + \cdots
\end{equation}

For the inner product, we define 
\begin{equation}
	g^A_{i,j} = <e_i,e_j>,
\end{equation}
and $g^{\gamma,\alpha}_A$ be its inverse.

Then, in coordinate, potential can be written as (with Einstein sign convention)
\begin{equation}
\Phi^A(\bx)=	\sum_{k \geq 1} \frac{1}{k+1} \,^Am^j_{i_1,\cdots,i_k} g_{j,i_{k+1}} x_{i_1}x_{i_2} \cdots x_{i_{k+1}}
\end{equation}

We define
 $h^*\Phi^A$ by the change of coordinate given by (\ref{changeco})
 
\begin{equation}
	h^*\Phi^A = \sum_{k \geq 1} \frac{1}{k+1} \,^Am^j_{i_1,\cdots,i_k} g_{j,i_{k+1}} h^{i_1}_{j_{11},\cdots,j_{1a_1}} \cdots
	h^{i_{k+1}}_{j_{(k+1)1},\cdots,j_{(k+1)a_{k+1}}} y_{j_{11}} \cdots y_{j_{(k+1)a_{k+1}}} 
\end{equation}

Here is the theorem due to Kajiura, in a translated form.
\begin{theorem}[Proposition 4.16 \cite{Kaj}]
Let $A$,$B$ be cyclic $\AI$-algebras and let $h: B \to A$ be a cyclic $\AI$-homomorphism.
Then, the potentials $\Phi^A,\Phi^B$ are related by
$$\Phi^B = (h^* \Phi^A)_c,$$
where $_c$ means cyclization.
\end{theorem} 
 
The proof in \cite{Kaj} was given by one line, which says that the proof follows from
the non-degeneracy of the symplectic structures. The proof suggested by
Kajiura is to check rather the identity $$(\Phi^B,\;) = ((h^*\Phi^A)_c,\;),$$
instead of the original one, where the bracket here is the poisson bracket from the symplectic structure. In writing down the detailed arguments following Kajiura's suggestion,
we found that signs work in a very subtle way and we will prove the theorem with a slightly different identity
to make the signs work.

Namely, for the proof only, we introduce a new variable $t_j$ corresponding to a variable $x_j$ whose degree  is
$N-2 -|x_j|'$, which in fact  equals the degree $|x_i|'$ when $g_{i,j} \neq 0$. Then consider
$$\overleftarrow{\frac{\partial}{\partial x_i} }g^{ij}t_j,$$
which has degree zero.

We first explain how this derivation works.
The basic rule is the Leibniz formula.
$$(x_1x_2)\overleftarrow{\frac{\partial}{\partial x_i} }g^{ij}t_j
= x_1\big( x_2 \overleftarrow{\frac{\partial}{\partial x_i} }g^{ij}t_j \big)+ \big( x_1 \overleftarrow{\frac{\partial}{\partial x_i} }g^{ij}t_j \big) x_2$$
$$=x_1\big( x_2 \overleftarrow{\frac{\partial}{\partial x_i} }g^{ij}t_j \big)+ (-1)^{|x_1||x_2|} x_2 \big( x_1 \overleftarrow{\frac{\partial}{\partial x_i} }g^{ij}t_j \big)$$

We may always put $t_j$ at the end by cyclically rotating the expression after the derivation.
For example, consider the following cyclic non-commutative polynomial 
$$F= \sum_{\sigma \in \ZZ / n \ZZ} \frac{a_{\sigma(1)\cdots \sigma(n)}}{n} x_1x_2\cdots x_n,$$
where cyclicity means the following equalities
$$a_{\sigma(i)\cdots \sigma(i-1)} = (-1)^{|x_{i-1}|(|x_{i}| + \cdots |x_{i-2}|)}
a_{\sigma(i-1)\cdots \sigma(i-2)}.$$
If all $x_i$'s are distinct, one can check that
\begin{equation}\label{pot1}
F \overleftarrow{\frac{\partial}{\partial x_i}}g^{ij}t_j 
= a_{\sigma(i)\cdots \sigma(i-1)} x_{i+1}x_{i+2}\cdots x_{i-1} \big( g^{ij}t_j \big)
\end{equation}

The case when some of the $x_i$'s are equal, also can be done by
Leibniz rule as in the ordinary calculus but except that after derivation, one can only cyclically rotate the expression.
The reason we need to be careful is due to the following derivation. Consider a cyclic polynomial
$$G = \frac{a}{n} (x_i)^n.$$

First, there may be two different notions of a cyclic polynomial, 
depending on whether one allows a monomial of type $(x_i)^n$ to be cyclic element always or to be so only if 
$$(x_i)^n = x_i(x_i)^{n-1} = (-1)^{(n-1)|x_i|^2} (x_i)^{n-1}x_i = (-1)^{(n-1)|x_i|^2} x_i^n,$$
which holds if $n$ is odd or $|x_i|$ is even. We call the former convention cyclic in a general sense, and the latter convention cyclic in a strict sense.
The potential defined above is strictly cyclic. But for a general non-commutative polynomial, its cyclization exists
only in the general sense because the above monomial does not have a cyclization in the strict sense.

In any case, if we assume $G$ to be strictly cyclic,
we have 
\begin{equation}\label{eq:der}
G  \overleftarrow{\frac{\partial}{\partial x_i}}g^{ij}t_j = a (x_i)^{n-1} g^{ij}t_j.
\end{equation}
Note that the results \ref{pot1}, \ref{eq:der} does not have the fraction $1/n$ anymore.
Without introducing the variable $t_j$, derivation would have resulted several different signs as
the derivative $\overleftarrow{\frac{\partial}{\partial x_i}}$ passes through $x_i$'s,
in which case we cannot get rid of the fraction $1/n$.
Only after we get rid of the fraction, we can apply various $\AI$-formulas.
This derivation may be thought as a way of cutting cyclic symmetry and putting a tag on the last element.
It is not hard to see  that the cyclic non-commutative polynomial with the given partial derivative is
determined uniquely up to a constant term.

Now we begin the proof of the Kajiura's theorem.
\begin{proof}
Instead of using $(\Phi^B,\;) = ((h^*\Phi^A)_c,\;)$, we will prove that
\begin{equation}\label{eq:inv}
\Phi^B \overleftarrow{\frac{\partial}{\partial y_\alpha}} g^{\alpha,\gamma}_B t_\gamma
= (h^*\Phi^A)_c \overleftarrow{\frac{\partial}{\partial y_\alpha}} g^{\alpha,\gamma}_B t_\gamma,
\end{equation}
which is enough due to the non-degeneracy of $g^{\alpha,\gamma}$.

First, consider the left hand side.
$$\Phi^B(y) = \sum_{k \geq 1} \frac{1}{k+1} \,^Bm^j_{i_1,\cdots,i_k} g^B_{j,i_{k+1}} y_{i_1}\cdots y_{i_{k+1}}$$
Hence, 
$$\Phi(y)^B \overleftarrow{\frac{\partial}{\partial y_\alpha}} g^{\alpha,\gamma}_B t_\gamma 
= \sum_{k} \,^Bm^j_{i_1,\cdots,i_k} g^B_{j,i_{k+1}} y_{i_1}\cdots y_{i_{k}} (y_{i_{k+1}} \overleftarrow{\frac{\partial}{\partial y_{i_{k+1}}}}) g^{i_{k+1},\gamma}_B t_\gamma $$
\begin{equation}\label{lhs}
 =  \sum_{k} \,^Bm^j_{i_1,\cdots,i_k}  y_{i_1}\cdots y_{i_{k}}  t_j.
\end{equation}
Note that at the first line, the fraction $1/k+1$ disappeared due to the derivation as explained above,
and the second equality follows by multiplying $g$ and $g^{-1}$.

Now, we consider the right hand side of the equation \ref{eq:inv}.
When we apply the derivative  
$\overleftarrow{\frac{\partial}{\partial y_\alpha}} g^{\alpha,\gamma}_B t_\gamma,$
we first organize the derivation into $(k+1)$ groups, corresponding to
the expressions $h^{i_1},\cdots,h^{i_{k+1}}$. Due to the cyclic symmetry of $\Phi^A$(in a strict sense),
the fraction $\frac{1}{k+1}$ is canceled by the derivation and  we have
$$(h^*\Phi^A)_c \overleftarrow{\frac{\partial}{\partial y_\alpha}} g^{\alpha,\gamma}_B y_\gamma
= \sum_{k}   \,^Am^j_{i_1,\cdots,i_k} g^A_{j,i_{k+1}}
 h^{i_1}_{j_{11},\cdots,j_{1a_1}} \cdots 	h^{i_{k+1}}_{j_{(k+1)1},\cdots,j_{(k+1)a_{k+1}}}$$
 
\begin{equation}\label{inv1}
	  \cdot y_{j_{11}} \cdots y_{j_{ka_{k}}} \big( (y_{j_{(k+1)1}} \cdots y_{j_{(k+1)a_{k+1}}}) 
 \overleftarrow{\frac{\partial}{\partial y_\alpha}} g^{\alpha,\gamma}_B t_\gamma, \big).
\end{equation}
 
Now, $\AI$-homomorphism relation may be written as
$$ \sum \,^Am^j_{i_1,\cdots,i_k} h^{i_1}_{j_{11},\cdots,j_{1a_1}} \cdots 	h^{i_{k}}_{j_{(k)1},\cdots,j_{(k)a_{k}}}$$
$$ = \sum h^j_{j_{11},\cdots,\delta,\cdots,j_{(k)a_{k}}} \,^Bm^\delta_{j_{*},\cdots,j_{**}}$$
Applying this formula to \ref{inv1},
we obtain
$$ =  \sum_{k} \big( h^j_{j_{11},\cdots,\delta,\cdots,j_{(k)a_{k}}} \,^Bm^\delta_{j_{*},\cdots,j_{**}} \big)
g^A_{j,i_{k+1}} h^{i_{k+1}}_{j_{(k+1)1},\cdots,j_{(k+1)a_{k+1}}}$$
\begin{equation}\label{inv2}
\cdot y_{j_{11}} \cdots y_{j_{ka_{k}}} \big( (y_{j_{(k+1)1}} \cdots y_{j_{(k+1)a_{k+1}}}) 
 \overleftarrow{\frac{\partial}{\partial y_\alpha}} g^{\alpha,\gamma}_B t_\gamma, \big).
\end{equation}

Now, the definition of a cyclic morphism (Lemma \ref{cycmordef}) implies that
the sum $$ \sum  h^j_{j_{11},\cdots,\delta,\cdots,j_{(k)a_{k}}}
g^A_{j,i_{k+1}} h^{i_{k+1}}_{j_{(k+1)1},\cdots,j_{(k+1)a_{k+1}}}$$
as the arguments of $h^j$ and $h^{i_{k+1}}$ varies over all possible two partitions of the set
$$( j_{11},\cdots,\delta,\cdots,j_{(k)a_{k}}, j_{(k+1)1},\cdots,j_{(k+1)a_{k+1}} ),$$
vanishes except the case when both $h^j, h^{i_{k+1}}$ have only one argument each.
The property $$<f_i,f_j>_B=<h(f_i),h(f_j)>_A$$ may be written as
$$ g^B_{i,j}=h^i_l g_{l,k}^A h^j_k$$
Hence, in our case, in the expression \ref{inv2} 
we have a relation $$g^B_{\delta,j_{(k+1)1}} = h^j_\delta g^A_{j,i_{k+1}} h^{i_{k+1}}_{j_{(k+1)1}}.$$
The expression \ref{inv2} becomes
$$\sum g^B_{\delta,j_{(k+1)1}} \,^Bm^\delta_{j_{n1},\cdots,j_{nl}} y_{j_{n1}} \cdots y_{j_{nl}}
\big( y_{j_{(k+1)1}}   \overleftarrow{\frac{\partial}{\partial y_\alpha}} g^{\alpha,\gamma}_B t_\gamma, \big)$$
$$ = \sum \,^Bm^\delta_{j_{n1},\cdots,j_{nl}} y_{j_{n1}} \cdots y_{j_{nl}} t_\delta $$
where we multiplied $g, g^{-1}$ to get the equality.
This agrees with the expression \ref{lhs}, and this proves the theorem.
\end{proof}

\section{Filtered cyclic $\AI$-algebra}
We introduce an analogous definition of cyclicity in the case of a filtered $\AI$-algebras,
which was introduced in the construction of the $\AI$-algebra of Lagrangian submanfolds by
Fukaya, Oh, Ohta and Ono \cite{FOOO}. 
 For the definition of filtered $\AI$-algebras, we refer readers to \cite{FOOO} for details.
Here we discuss filtered $\AI$-algebras over the filtered ring $\NOV$, which is the universal Novikov ring.
The universal Novikov ring is defined as ($T,e$ as formal parameters)
$$\NOV = \{ \sum_{i=0}^\infty a_i T^{\lambda_i} e^{n_i}| a_i \in \kk, \lambda_i \in \RR, n_i \in \ZZ \,\; \textrm{and} \;\; \lim_{i \to \infty} \lambda_i = \infty 
\}.$$
We have a graded subring $\Lambda_{0,nov}$ if all $\lambda_i \geq 0$, and the filtration $F^{\lambda}$ is defined by
the condition $\lambda_i \geq \lambda$.
Consider the filtered $\AI$-algebra structure on $(C,m_{\geq 0})$.
In the filtered $\AI$-algebras, there may exists $m_0$ term but is assume to satisfy
$m_0(1) \in F^{\lambda}C$ for some $\lambda >0$.
And, $\overline{C}$ is a $\kk$-vector space with $\overline{C} \otimes_{\kk}\Lambda_{0,nov} \cong C$.

 We consider only the gapped filtered $\AI$-algebras, where the gapped condition
is defined as follows.
The monoid $G \subset \RR_{\geq 0} \times 2 \ZZ$ is assumed to satisfy the following conditions
\begin{enumerate}
\item The projection $\pi_1(G) \subset \RR_{\geq 0}$ is discrete.
\item $G \cap (\{0\} \times 2\ZZ) = \{(0,0)\}$
\item $G \cap (\{\lambda \} \times 2\ZZ)$ is a finite set for any $\lambda$.
\end{enumerate}
Then, $(C,m_{\geq 0})$ is said to be $G$-gapped if there exists $\kk$-vector space
homomorphisms $m_{k,\beta}:B_k \overline{C}[1] \to \overline{C}[1]$ for $k=0,1,\cdots,$ $\beta=(\lambda(\beta),\mu(\beta)) \in G$
such that
$$m_k = \sum_{\beta \in G} T^{\lambda(\beta)}e^{\mu(\beta)/2} m_{k,\beta}.$$

Let $(C,\{m_{i,\beta} \}_{i \geq 0, \beta \in G})$ be a gapped filtered $\AI$-algebra as above. 

\begin{definition}
A filtered gapped $\AI$-algebra $(C,\{m_*\})$ is said to have a {\it cyclic symmetric} inner product if
there exists a skew-symmetric non-degenerate, bilinear map $$<,> : \overline{C}[1] \otimes \overline{C}[1] \to \kk,$$
which is extended linearly over $C$, 
such that for all integer $k \geq 0$, $\beta \in G$, 
\begin{equation}
	<m_{k,\beta}(x_1,\cdots,x_k),x_{k+1}> = (-1)^{K}<m_{k,\beta}(x_2,\cdots,x_{k+1}),x_{1}>.
\end{equation}
 where $K = |x_1|'(|x_2|' + \cdots +|x_{k+1}|')$.
\end{definition}

The definitions, and theorems in this paper for unfiltered $\AI$-algebras can be extended to those of gapped filtered $\AI$-algebras.
\begin{definition}
An $\AI$-homomorphism $\{h_k\}_{k\geq 1 }$ between two gapped(over G) filtered $\AI$-algebras with cyclic symmetric inner products is called
a cyclic gapped filtered $\AI$-homomorphism if 
\begin{enumerate}
\item $h_{1,0}$ preserves inner product $<a,b> = <h_{1,0}(a),h_{1,0}(b)>$.
\item \begin{equation}
\sum_{i+j=k, \beta_1+\beta_2 =\beta \in G} <h_{i,\beta_1}(x_1,\cdots,x_i), h_{j,\beta_2}(x_{i+1},\cdots,x_k)> =0.	
\end{equation}
\end{enumerate}
\end{definition}
\begin{remark}
This is equivalent to the definition using the commutative diagram as before.
\end{remark}

Here $h_{1,0} = h_{1,\beta_0}$ where $\beta_0$ is the zero in the monoid $G$, or the homotopy class of
the constant holomorphic maps in the case of \cite{FOOO}.

\begin{corollary}
The characterization theorems \ref{thm:hi}, \ref{thm:sh} also holds for the cyclic gapped filtered $\AI$-algebra
($\AI$-bimodules).
\end{corollary}
\begin{proof}
The proof has been written so that the analogous arguments in the gapped filtered case also holds.
The difference is on the induction. In the gapped filtered case, the induction
should be run over the sum of two indices. As $\lambda(\beta)$ is
discrete, we can find an increasing sequence $\lambda_i$ with $\lim \lambda_i  = \infty$ which covers
the image of $\lambda(G) \subset \RR_{\geq 0}$.
We run the induction over the sum $k+i$, where for the inductive step $k+i=N$,
we consider the terms $m_{k,\beta}$ with $\lambda(\beta) = \lambda_i$.
\end{proof}

\begin{theorem}\label{thm:pofil}
Let $A$,$B$ be cyclic gapped filtered $\AI$-algebras and let $h: B \to A$ be a cyclic gapped filtered $\AI$-homomorphism.
Then, the potentials $\Phi^A,\Phi^B$ are related by
$$\Phi^B = (h^* \Phi^A)_c + \;\textrm{constant},$$
where $_c$ means cyclization.
\end{theorem} 
\begin{proof}
The proof given in the previous section can be easily generalized to
the filtered case except that in the filtered case, the related equivalence is up to a constant addition.
This is because, the change of parameter may produce constant terms (comming from $h_0$). But when we prove an identity,
we check the identity after taking a suitable derivation. Hence the proof does not show what happens for the constant term.
But in general, this constant term is expected to change when we change over $\AI$-homomorphism. 
We do not know if there is a way to express and preserve the constant contribution from the potential, which we leave as an
open question.
\end{proof}
\bibliographystyle{amsalpha}

\end{document}